\documentclass[conference, 10pt, twocolumn]{ieeeconf} 
\pdfoutput=1

\def\BibTeX{{\rm B\kern-.05em{\sc i\kern-.025em b}\kern-.08em
    T\kern-.1667em\lower.7ex\hbox{E}\kern-.125emX}}
    
\usepackage[left=54pt, right=54pt,bottom=54pt, top=54pt]{geometry}

\usepackage{amsmath,mathrsfs,amsfonts,amssymb,graphicx,epsfig}
 \usepackage{amsthm}
\usepackage{subcaption}
\usepackage{color,multirow,rotating}
\usepackage{algorithm,algpseudocode,algorithmicx}
\usepackage{cite,url,framed,bm,balance}

\newtheorem{theorem}{Theorem}

\newtheorem{proposition}{Proposition}
\newtheorem{remark}{Remark}

\usepackage{tikz}
\usetikzlibrary{calc,trees,positioning,arrows,chains,shapes.geometric,decorations.pathreplacing,decorations.pathmorphing,shapes, matrix,shapes.symbols}

\newcommand{\differential}{{\rm{d}}}

\newcommand{\blkdiag}{\mathrm{bldiag}}

\makeatletter
\newcommand{\dotminus}{\mathbin{\text{\@dotminus}}}

\newcommand{\@dotminus}{%
  \ooalign{\hidewidth\raise1ex\hbox{.}\hidewidth\cr$\m@th-$\cr}%
}

\allowdisplaybreaks

\title{\LARGE\textbf{
Certifying the Intersection of Reach Sets of Integrator Agents\\with Set-valued Input Uncertainties}
}

\author{Shadi Haddad and Abhishek Halder
\thanks{Shadi Haddad and Abhishek Halder are with the Department of Applied Mathematics, University of California, Santa Cruz, CA 95064, USA, {\tt\small{\{shhaddad,halder\}@ucsc.edu}}.%
}}

\IEEEoverridecommandlockouts
\begin{document}

\maketitle
\pagenumbering{arabic}

\begin{abstract}
We consider the problem of verifying safety for a pair of identical integrator agents in continuous time with compact set-valued input uncertainties. We encode this verification problem as that of certifying or falsifying the intersection of their reach sets. We transcribe the same into a variational problem, namely that of minimizing the support function of the difference of the two reach sets over the unit sphere. We illustrate the computational tractability of the proposed formulation by developing two cases in detail, viz. when the inputs have time-varying norm-bounded and generic hyperrectangular uncertainties. We show that the latter case allows distributed certification via second order cone programming.
\end{abstract}


\section{Introduction}\label{sec:introduction}
We consider a pair of integrator agents labeled as \texttt{A} and \texttt{B}, each with $n$ states, $m$ inputs and relative degree $\bm{r}=(r_1,\hdots,r_m)^{\top}\in\mathbb{N}^{m}$ where $r_{1}+\hdots + r_{m} = n$. For example, when $n=5$, $m=2$, and $\bm{r}=(3,2)^{\top}$, then the individual integrator agent dynamics in continuous time is of the form
\begin{align}
\dot{\bm{x}}^{i} = \left(\begin{array}{@{}ccc|cc@{}}
0 & 1 & 0 & 0 & 0\\
0 & 0 & 1 & 0 & 0\\
0 & 0 & 0 & 0 & 0\\ \hline
0 & 0 & 0 & 0 & 1\\
0 & 0 & 0 & 0 & 0\\
\end{array}\right)
\bm{x}^{i} + \left(\begin{array}{@{}c|c@{}}
0 & 0\\
0 & 0\\ 
1 & 0\\\hline
0 & 0\\
0 & 1
\end{array}\right)\bm{u}^{i},    
\label{Example3block2block}
\end{align}
wherein for $i\in\{\texttt{A},\texttt{B}\}$, the vectors $\bm{x}^{i},\bm{u}^{i}$ respectively denote the state and input of the agents \texttt{A} and \texttt{B}. 

In general, for $i\in\{\texttt{A},\texttt{B}\}$, the integrator agent dynamics takes the form
\begin{align}
\dot{\bm{x}}^{i} \!= \!\underbrace{\blkdiag\left(\bm{A}_{1},\hdots,\bm{A}_{m}\right)}_{=:\bm{A}}\!\bm{x}^{i} \!+\! \underbrace{\blkdiag\left(\bm{b}_{1},\hdots,\bm{b}_{m}\right)}_{=:\bm{B}}\!\bm{u}^{i},    
\label{IntegratorODE}    
\end{align}
where $\blkdiag\left(\cdot\right)$ denotes block diagonal matrix whose arguments constitute its diagonal blocks, and
\begin{align}
\boldsymbol{A}_{j}:=\left(\mathbf{0}_{r_{j} \times 1}\left|\boldsymbol{e}_{1}^{r_{j}}\right| \boldsymbol{e}_{2}^{r_{j}}|\ldots| \boldsymbol{e}_{r_{j}-1}^{r_{j}}\right), \quad \boldsymbol{b}_{j}:=\boldsymbol{e}_{r_{j}}^{r_{j}},
\label{blkdiagAB}    
\end{align}
for $j \in [m]:=\{1,2,\hdots,m\}$. In (\ref{blkdiagAB}), the notation $\mathbf{0}_{r_{j} \times 1}$ stands for the vector of zeros of size $r_{j} \times 1$. For $k\leq\ell$, we use $\bm{e}_{k}^{\ell}$ to denote the $k$th standard basis column vector in $\mathbb{R}^{\ell}$.

In this work, we consider the problem of checking whether the agents \texttt{A} and \texttt{B}, starting from respective initial conditions $\bm{x}_{0}^{\texttt{A}},\bm{x}_{0}^{\texttt{B}}\in\mathbb{R}^{n}$, and respective input uncertainties modeled as compact\footnote{This means that the set-valued trajectories in \eqref{InputSets} are continuous w.r.t. $s$, and the sets are compact for all $0\leq s \leq t$.} input set-valued trajectories
\begin{align}
\mathcal{U}^{\texttt{A}}(s),\,\mathcal{U}^{\texttt{B}}(s)\subset\mathbb{R}^{m}, \quad 0\leq s\leq t,
\label{InputSets}
\end{align}
may result in \emph{intersecting} reach sets $\mathcal{X}_{t}^{\texttt{A}}, \mathcal{X}_{t}^{\texttt{B}}\subset\mathbb{R}^{n}$ at a given time $t$. Formally, the (forward) reach sets are given by
\begin{align}
&\!\mathcal{X}_{t}^{i}:=\!\!\!\bigcup_{\text{measurable}\;\bm{u}^{i}(\cdot)\in\mathcal{U}^{i}(\cdot)}\!\{\bm{x}^{i}(t)\in\mathbb{R}^{n} \mid \eqref{IntegratorODE}, \; \bm{x}^{i}(t=0)=\bm{x}_{0}^{i}, \nonumber\\
&\bm{u}^{i}(s)\in\mathcal{U}^{i}(s)\;\text{compact for all}\; 0\leq s\leq t\},\; i\in\{\texttt{A},\texttt{B}\},
\label{DefReachSets}    
\end{align}
such that $\bm{x}_{0}^i=(\bm{x}_{10}^i,\hdots,\bm{x}_{m0}^i)^{\top}\subset\mathbb{R}^{n}$. In words, the forward reach set at time $t$ is the set of states that the individual agent may reach subject to its integrator dynamics, given initial condition and compact set-valued input uncertainties. In this paper, we only consider reach sets forward in time, and hereafter refer to the same as reach sets.

It is well-known \cite{varaiya2000reach} that the reach sets $\mathcal{X}_{t}^{i}$ in \eqref{DefReachSets} are \emph{compact convex} provided the input sets \eqref{InputSets} are \emph{compact}. Furthermore, the sets $\mathcal{X}_{t}^{i}$ are invariant under the closure of convexification of the input sets \eqref{InputSets}; see e.g., \cite[Prop. 6.1, Thm. 6.3]{yong1999stochastic}. We would like to certify (or falsify) if $\mathcal{X}_{t}^{\texttt{A}}\cap\mathcal{X}_{t}^{\texttt{B}} \neq (=) \varnothing$ at a given time $t$.

Collision detection and avoidance for integrator agents have received attention \cite{cao2010distributed,abdessameud2013consensus,li2011finite,ajorlou2012bounded} in the robotics literature. Beyond serving as simple models for agent dynamics, integrators also arise as Brunovsky normal forms of feedback linearizable systems resulting from a diffeomorphic change of state coordinates. On the other hand, detecting intersection of reach sets amounts to verifying safe operation in the sense of checking separability of the agents' states or lack thereof subject to set-valued uncertainties. In particular, detecting the intersection of reach sets of \emph{static} state feedback linearizable systems, can be shown\footnote{See Supplementary Material} to be equivalent to detecting the same in the normal coordinates. With this motivation, the present work investigates certifying the intersection of reach sets of integrator agents using ideas from convex geometry. 

The organization of this paper is as follows. In Sec. \ref{sec:SptFnPrelim}, we provide some background on the support function of a compact convex set in general, and that of the integrator reach set in particular. In Sec. \ref{sec:formulation}, we formulate the problem of detecting intersection among the reach sets of two integrator agents, labeled as \texttt{A} and \texttt{B}, as that of solving a nonconvex problem involving their support functions. We next show that when the input uncertainty sets are time-varying norm balls (Sec. \ref{sec:NormBallInputSets}) or time-varying generic hyperrectangles (Sec. \ref{sec:HyperrectangleInputSets}), then a lossless convexification of the nonconvex formulation solves the intersection detection problem. 

\section{Background on Support Function}\label{sec:SptFnPrelim}
 The support function $h_{\mathcal{K}}(\cdot)$ of a compact convex set $\mathcal{K} \subset \mathbb{R}^{n}$, is given by
\begin{align}
h_{\mathcal{K}}(\bm{y}) := \underset{\bm{x}\in\mathcal{K}}{\sup}\:\{\langle\bm{y},\bm{x}\rangle \mid \bm{y}\in\mathbb{S}^{n-1}\},
\label{DefSptFn}	
\end{align}
where $\langle\cdot,\cdot\rangle$ denotes the standard Euclidean inner product, and $\mathbb{S}^{n-1}$ is the unit sphere embedded in $\mathbb{R}^{n}$. 
For details on the support function, we refer the readers to \cite[Ch. V]{hiriart2013convex}.

%


The following result (proof in Appendix \ref{AppendixProofPro:proSptFn}) will come in handy in the ensuing development.

\begin{proposition}\label{proSptFn}
Consider the integrator dynamics with $n$ states, $m$ inputs, and relative degree vector $\bm{r}=(r_1,\hdots,r_m)^{\top}\in\mathbb{N}^{m}$. For compact input uncertainty sets $\mathcal{U}(s)\subset\mathbb{R}^{m}$, $0\leq s \leq t$, the support function of the integrator reach set $\mathcal{X}_t$ at time $t$ starting from the initial condition $\bm{x}_{0} = \left(\bm{x}_{10},\bm{x}_{20},\hdots,\bm{x}_{m0}\right)^{\top}\in \mathbb{R}^n$, is given by
\begin{align}
\label{SptFnIntegratorpGen}
h_{\mathcal{X}_t}\left(\bm{y}\right) &=
\sum_{j=1}^{m}~ \langle\bm{y}_{j},\exp\left(t\bm{A}\right)\bm{x}_{j0}\rangle\nonumber\\ 
&\quad+ \int^t_0 \underset{\bm{u}(s)\in\mathcal{U}(s) }{\sup} \sum_{j=1}^{m}\langle\bm{y}_{j},\bm{\xi}_{j}(s)\rangle\: u_{j}(s)\: \differential{s},
\end{align}
where the subvector $\bm{x}_{j0}\in\mathbb{R}^{r_j}$, $j\in[m]$, and
\begin{align}
\bm{\xi}_{j}(s) := \!\begin{pmatrix}
s^{r_j-1}/(r_j-1)! \\
s^{r_j-2}/(r_j-2)! \\
  \vdots\\
  s\\
   1
\end{pmatrix},
~\bm{\xi}(s):=\!\begin{pmatrix}
\bm{\xi}_1(s)\\
\vdots\\
\bm{\xi}_m(s)
\end{pmatrix}.  
\label{xiVector}
\end{align}
\end{proposition}
The proof in Appendix \ref{AppendixProofPro:proSptFn} reveals a relationship between the support function of the reach set $\mathcal{X}_{t}$ with that of the compact input sets $\mathcal{U}(s)$, $0\leq s \leq t$, given by
\begin{align}
h_{\mathcal{X}_t}(\bm{y}) = \langle\bm{y},\exp(t\bm{A})\bm{x}_{0}\rangle \!+ \!\!\int_{0}^{t}\!\!\!h_{\mathcal{U}(s)}\!\!\left(\!\left(\exp(s\bm{A})\bm{B}\right)^{\top}\!\!\bm{y}\!\right)\!\differential s.	
\end{align}
When $\mathcal{U}(s)$ is compact nonconvex, then $h_{\mathcal{U}(s)}(\cdot)$ is the support function of the closure of the convex hull of $\mathcal{U}(s)$. 

\section{Problem Formulation}\label{sec:formulation}
One way to check if $\mathcal{X}_{t}^{\texttt{A}}\cap\mathcal{X}_{t}^{\texttt{B}} = \varnothing$ or not, is to compute the distance between the reach sets $\mathcal{X}_{t}^{\texttt{A}},\mathcal{X}_{t}^{\texttt{B}}$, given by
\begin{align}
\label{DistanceBetweenSets} 
{\rm{dist}}\left(\texttt{A},\texttt{B}\right) := \underset{\stackrel{}{\bm{x}^{\texttt{A}}\in \mathcal{X}_{t}^{\texttt{A}},\bm{x}^{\texttt{B}}\in\mathcal{X}_{t}^{\texttt{B}}}}{\min}\;\|\bm{x}^{\texttt{A}}-\bm{x}^{\texttt{B}}\|_{2}^{2}.   
\end{align}
Clearly, ${\rm{dist}}\left(\texttt{A},\texttt{B}\right) > 0$ if and only if $\mathcal{X}_{t}^{\texttt{A}}\cap\mathcal{X}_{t}^{\texttt{B}} = \varnothing$, and $=0$ otherwise. However, computing the distance between the reach sets $\mathcal{X}_{t}^{\texttt{A}},\mathcal{X}_{t}^{\texttt{B}}$ requires analytic handle on the boundary of these reach sets, which can be computationally difficult depending on the geometry of the input sets \eqref{InputSets}.

To circumvent this difficulty, we take an alternative approach 
based on the difference set of $\mathcal{X}_{t}^{\texttt{A}}$ and $\mathcal{X}_{t}^{\texttt{B}}$, given by the compact convex set  
\begin{align}
\mathcal{X}_{t}^{\texttt{A}}\dotminus\mathcal{X}_{t}^{\texttt{B}} := \{\bm{x}^{\texttt{A}}-\bm{x}^{\texttt{B}} \mid \bm{x}^{\texttt{A}}\in \mathcal{X}_{t}^{\texttt{A}},\bm{x}^{\texttt{B}}\in\mathcal{X}_{t}^{\texttt{B}}\}.    
\label{DefMinkowskiDifference}    
\end{align}
Notice that $\mathcal{X}_{t}^{\texttt{A}}\dotminus\mathcal{X}_{t}^{\texttt{B}} = \mathcal{X}_{t}^{\texttt{A}}\dotplus\left(-\mathcal{X}_{t}^{\texttt{B}}\right)$ where $\dotplus$ denotes the Minkowski sum, and $\mathcal{X}_{t}^{\texttt{A}}\dotminus\mathcal{X}_{t}^{\texttt{B}}$ is not the same as the Minkowski a.k.a. Pontryagin difference \cite[p. 139]{schneider2014convex}, \cite{montejano1996some}.

Checking the intersection between $\mathcal{X}_{t}^{\texttt{A}}$ and $\mathcal{X}_{t}^{\texttt{B}}$, is then equivalent to verifying if the zero vector $\bm{0}\in\mathbb{R}^{n}$ belongs to the set \eqref{DefMinkowskiDifference}. This can in turn be related \cite{zheng2015generalized,hornus2017detecting} to conditions on the support function $h_{\mathcal{X}_{t}^{\texttt{A}}\dotminus\mathcal{X}_{t}^{\texttt{B}}}(\cdot)$, because
\begin{subequations}
\begin{align}
\!\!\!\!\!\!&\mathcal{X}_{t}^{\texttt{A}}\cap\mathcal{X}_{t}^{\texttt{B}}\neq\varnothing\:\Leftrightarrow\:\bm{0}\!\in\! \{\mathcal{X}_{t}^{\texttt{A}}\dotminus\mathcal{X}_{t}^{\texttt{B}}\}\nonumber\\ 
&\hspace*{0.85in}\Leftrightarrow\: \forall\;\bm{y}\in\mathbb{S}^{n-1}, h_{\mathcal{X}_{t}^{\texttt{A}}\dotminus\mathcal{X}_{t}^{\texttt{B}}}(\bm{y}) \geq 0, \label{MinkDiffIntersect}\\   
\!\!\!\!\!\!&\mathcal{X}_{t}^{\texttt{A}}\cap\mathcal{X}_{t}^{\texttt{B}}=\varnothing\:\Leftrightarrow\: \bm{0}\!\notin\! \{\mathcal{X}_{t}^{\texttt{A}}\dotminus\mathcal{X}_{t}^{\texttt{B}}\}\nonumber\\ 
&\hspace*{0.85in}\Leftrightarrow\: \exists\;\bm{y}\in\mathbb{S}^{n-1}\,\text{such that}\: h_{\mathcal{X}_{t}^{\texttt{A}}\dotminus\mathcal{X}_{t}^{\texttt{B}}}(\bm{y}) < 0. \label{MinkDiffDisjoint} 
\end{align}
\label{MinkowskiDiffIntersectionOrNot} 
\end{subequations}
Thus motivated, we propose certifying or falsifying the reach set intersection by computing
\begin{align}
\label{SptFnMinkDiffOptimization}  
\underset{\bm{y}\in\mathbb{S}^{n-1}}{\min}\: h_{\mathcal{X}_{t}^{\texttt{A}}\dotminus\mathcal{X}_{t}^{\texttt{B}}}(\bm{y}).     
\end{align}
Specifically, \eqref{SptFnMinkDiffOptimization} $\geq (<) 0 \:\Leftrightarrow\: \mathcal{X}_{t}^{\texttt{A}}\cap\mathcal{X}_{t}^{\texttt{B}} \neq (=) \varnothing$. In other words, the optimal cost \eqref{SptFnMinkDiffOptimization} can be used to certify or falsify the reach set intersection.

\begin{remark}\label{remark:supportfnoptimization}
In general, which of the two problems \eqref{DistanceBetweenSets} and \eqref{SptFnMinkDiffOptimization} is computationally more tractable, depends on the sets $\mathcal{X}_{t}^{\texttt{A}},\mathcal{X}_{t}^{\texttt{B}}$. For instance, when $\mathcal{X}_{t}^{\texttt{A}},\mathcal{X}_{t}^{\texttt{B}}$ are ellipsoids, simple algorithms are known \cite{lin2002distance} for solving \eqref{DistanceBetweenSets}, but solving \eqref{SptFnMinkDiffOptimization} requires heavier computation \cite{iwata2015computing}. In the ensuing sections, we explain how detecting the intersection between the integrator reach sets $\mathcal{X}_{t}^{\texttt{A}},\mathcal{X}_{t}^{\texttt{B}}$ using \eqref{SptFnMinkDiffOptimization} turns out to be computationally more benign than \eqref{DistanceBetweenSets}.
\end{remark}

Since support function is distributive over the Minkowski sum, we have \begin{align}
h_{\mathcal{X}_{t}^{\texttt{A}}\dotminus\mathcal{X}_{t}^{\texttt{B}}}(\bm{y}) = h_{\mathcal{X}_{t}^{\texttt{A}}}(\bm{y}) + h_{-\mathcal{X}_{t}^{\texttt{B}}}(\bm{y}).
\label{Intermed}    
\end{align}
From the definition of support function, we also have
\begin{align}
\label{SptFnOfNegation}
\!\!\!\!h_{-\mathcal{X}_{t}^{\texttt{B}}}(\bm{y}) =\! \underset{\bm{x}\in-\mathcal{X}_{t}^{\texttt{B}}}{\sup}\langle\bm{y},\bm{x}\rangle &= \underset{\bm{x}\in\mathcal{X}_{t}^{\texttt{B}}}{\sup}\langle\bm{y},-\bm{x}\rangle \nonumber\\  &=\underset{\bm{x}\in\mathcal{X}_{t}^{\texttt{B}}}{\sup}\langle-\bm{y},\bm{x}\rangle = h_{\mathcal{X}_{t}^{\texttt{B}}}(-\bm{y}). 
\end{align}
Using \eqref{Intermed} and \eqref{SptFnOfNegation}, we rewrite \eqref{SptFnMinkDiffOptimization} as
\begin{align}
\label{SptFnOptimizationCollisionDetection} 
\underset{\|\bm{y}\|_{2}=1}{\min}\: h_{\mathcal{X}_{t}^{\texttt{A}}}(\bm{y}) + h_{\mathcal{X}_{t}^{\texttt{B}}}(-\bm{y}).     
\end{align}
Recall that a support function is convex in its argument, and convex function composed with an affine map remains convex. Thus, the objective in \eqref{SptFnOptimizationCollisionDetection} is a sum of convex functions, and hence convex. Because $\mathbb{S}^{n-1}$ is compact, by Weirstrass extreme value theorem, \eqref{SptFnOptimizationCollisionDetection} admits global minimum. Checking whether $\mathcal{X}_{t}^{\texttt{A}}\cap\mathcal{X}_{t}^{\texttt{B}} = \varnothing$ or not, reduces to checking the sign of the minimum in \eqref{SptFnOptimizationCollisionDetection}.

We next develop these ideas for two specific choices of time-varying set-valued input uncertainties: norm bounded (Sec. \ref{sec:NormBallInputSets}) and hyperrectangular uncertainty sets (Sec. \ref{sec:HyperrectangleInputSets}).

\section{The Case when $\mathcal{U}^{\texttt{A}},\mathcal{U}^{\texttt{B}}$  are Norm Balls}\label{sec:NormBallInputSets}
In this section, we consider the time-varying norm-bounded input uncertainty sets
\begin{align}
\label{inputNorm}
\mathcal{U}(s):= \Big\{ \bm{u}(s) \in \mathbb{R}^m ~\big\vert~ \| \bm{u}(s) \|_p \leq \ell(s) \Big\}\:\text{for all}~s\in [0,t],    
\end{align}
where $\|.\|_{p}$ denotes the $p$-norm for $1\leq p \leq \infty$, and $\ell:[0,t]\mapsto\mathbb{R}_{>0}$ is a given smooth function. In \eqref{inputNorm}, without loss of generality, we consider the norm ball to be symmetric about the origin, since a translation of the norm ball does not change the arguments provided next. By specializing \eqref{SptFnIntegratorpGen}, we get the following result (proof in Appendix \ref{AppendixProofThm:Normspf}).
\begin{theorem}\label{thm:Normspf}
Consider the integrator dynamics with $n$ states, and $m$ inputs. For compact input uncertainty sets $\mathcal{U}(s)\subset\mathbb{R}^{m}$ given by \eqref{inputNorm} for all $0\leq s \leq t$, and initial condition $\bm{x}_0 \in \mathbb{R}^n$, the support function of the integrator reach set $\mathcal{X}_t$ at time $t$ is
\begin{align}
h_{\mathcal{X}_{t}}\left(\bm{y}\right) &= 
\sum_{j=1}^{m}~ \langle\bm{y}_{j},\exp\left(t\bm{A}\right)\bm{x}_{j0}\rangle\nonumber\\
&\quad+ \int^t_0 \ell(s)\|(\exp(s\bm{A})\bm{B})^{\top}\bm{y} \|_q  \:\differential{s},
\label{SptFnIntegratorpNorm}
\end{align}	
where $q$ is the H\"{o}lder conjugate of $p$, i.e., $\frac{1}{p}+\frac{1}{q}=1$.
\end{theorem}
\begin{remark}
We mentioned in Sec. \ref{sec:introduction} that the reach set $\mathcal{X}_t$ resulting from compact $\mathcal{U}(s)$ is the same as that resulting from the closure of the convex hull of $\mathcal{U}(s)$, $0\leq s\leq t$. Consequently, if the $p$ in \eqref{inputNorm} satisfies $0<p<1$, thus making the input norm balls nonconvex, then the corresponding reach sets will coincide with that resulting from the $p=1$ norm ball input uncertainty sets. This allows the effective domain of $p$ in \eqref{inputNorm} to be $(0,\infty]$.
\end{remark}

We next detail how having analytic handle on the support function of the reach set as in \eqref{SptFnIntegratorpNorm} can help detect reach set intersection among integrator agents using \eqref{SptFnOptimizationCollisionDetection}.

\subsection{Lossless Convexification}\label{subsec:convexification}
We suppose that the integrator agents \texttt{A} and \texttt{B} have input uncertainty sets $\mathcal{U}^{\texttt{A}}(s),\mathcal{U}^{\texttt{B}}(s)$ as in \eqref{inputNorm} with same $p$, respective bounds ${\ell}^{\texttt{A}}(s),{\ell}^{\texttt{B}}(s)$, and respective initial conditions $\bm{x}_{0}^{\texttt{A}},\bm{x}_{0}^{\texttt{B}}\in\mathbb{R}^{n}$.

The associated problem \eqref{SptFnOptimizationCollisionDetection} is nonconvex due to the unit sphere constraint $\|\bm{y}\|_{2}=1$. We convexify the same by relaxing it to the unit ball constraint $\|\bm{y}\|_{2}\leq 1$. Since ${\ell}^{\texttt{A}}(s),{\ell}^{\texttt{B}}(s)$ are positive for all $0\leq s \leq t$, the convexified version of \eqref{SptFnOptimizationCollisionDetection} becomes
\begin{align}\vspace*{-0.2in}
\!\underset{\bm{y}\in\mathbb{R}^{n},\|\bm{y}\|_{2}\leq 1}{\min}\: \underbrace{\langle\bm{c}(t),\bm{y}\rangle + \!\!\int_{0}^{t}\!\!\|\left(\bm{G}(s)\right)^{\top}\bm{y}\|_q \:\differential s}_{=:f_{0}(\bm{y})}, 
\label{ConvexifiedSubproblemNorm}    
\end{align}
where \vspace*{-0.2in}
\begin{subequations}
\begin{align}
\bm{c}(t) &:= \exp(t\bm{A})\!\left(\bm{x}_{0}^{\texttt{A}} - \bm{x}_{0}^{\texttt{B}}\right),\label{DefcNorm}\\
\bm{G}(s) &:= \left(\ell^{\texttt{A}}(s)+\ell^{\texttt{B}}(s)\right)\exp(s\bm{A})\bm{B}.\label{DefgammaNorm}
\end{align}
\label{DefcANDgammaNorm}
\end{subequations}
We approximate the integral in \eqref{ConvexifiedSubproblemNorm} w.r.t. $s$ via trapezoidal approximation\footnote{The uniform trapezoidal approximation with local truncation error $\mathcal{O}((\Delta s)^{3})$ may be replaced by other approximations such as the three point Simpson's rule with local truncation error $\mathcal{O}((\Delta s)^{4})$. 
While the accuracy of different numerical approximations for the integral may vary depending on the choice of approximation but the nature of the resulting optimization problems will remain the same.} with uniform step-size $\Delta s > 0$. In particular, uniformly discretizing $[0,t]$ into $K\in\mathbb{N}$ intervals with breakpoints $s_{k} = k\Delta s$ for $k=0,1,\hdots,K$, where $\Delta s := t/K$, results in the trapezoidal approximation
\begin{align}
 &\int_{0}^{t}\!\!\!\|\left(\bm{G}(s)\right)^{\top}\bm{y}\|_{q}\:\differential s \approx\nonumber\\
 &\qquad\frac{\Delta s}{2}\!\sum_{k=1}^{K}\!\left(\|\left(\bm{G}(s_{k-1})\right)^{\top}\bm{y} \|_q \!+\! \|\left(\bm{G}(s_k)\right)^{\top}\bm{y} \|_q\right).
 \label{TrapzApprox}
\end{align}
Letting 
\begin{align*}
\|\left(\bm{G}(s_{k})\right)^{\top}\bm{y}\|_{q}\leq \theta_{k}, \quad k=0,1,\hdots,K,
\end{align*}
and $\bm{\theta}:=(\theta_0,\theta_1,\hdots,\theta_{K})\in\mathbb{R}^{K+1}_{\geq 0}$, we next define
{\small{\begin{subequations}
\begin{align}
&\bm{\eta} := \begin{pmatrix}
\bm{y}\\
\bm{\theta}
\end{pmatrix}\in\mathbb{R}^{n+K+1},\label{defzjNorm}\\
&\bm{\omega} := \Delta s\begin{pmatrix}
1/2\\
\bm{1}_{K-1}\\
1/2
\end{pmatrix}\in\mathbb{R}^{K+1}_{>0}
,\label{defomegajNorm}\\ 
&\bm{\kappa}(t) := \begin{pmatrix}
\bm{c}(t)\\
\bm{\omega}
\end{pmatrix}\in\mathbb{R}^{n+K+1}, \label{defelljNorm}\\
&\vspace*{-0.8 in}\bm{M}_k := \begin{pmatrix}
\begin{array}{c|c}
\bm{G}(s_k)^{\top} & \bm{0}_{m\times(K+1)} \\
\end{array} 
\end{pmatrix}\in\mathbb{R}^{m\times(n+K+1)}\\
&\bm{N} := \!\begin{pmatrix}
\begin{array}{c|c}
\bm{I}_{n} & \bm{0}_{n\times(K+1)}
\end{array} 
\end{pmatrix}\in\mathbb{R}^{n\times(n+K+1)}.\label{defNjNorm}\\
&\tilde{\bm{N}} := \!\begin{pmatrix}
\begin{array}{c|c}
\bm{0}_{(K+1)\times n} & \bm{I}_{K+1}
\end{array} 
\end{pmatrix}\in\mathbb{R}^{(K+1)\times(n+K+1)}.\label{deftildeNjNorm} 
\end{align}
\label{defSOCPparamNorm}
\end{subequations}}}
In \eqref{defSOCPparamNorm}, the symbols $\bm{1}$, $\bm{0}$, $\bm{I}$ respectively denote the array of ones, zeros and identity matrix of appropriate sizes.

With the above variable definitions in hand, we transcribe \eqref{ConvexifiedSubproblemNorm} into the epigraph form
\begin{subequations}
\begin{align}
&\underset{\bm{\eta}\in\mathbb{R}^{n+K+1}}{\min}\:\langle \bm{\kappa}(t),\bm{\eta}\rangle   \label{SOCPobjNorm}\\ 
&\text{subject to} \nonumber\\
&\|\bm{M}_k\bm{\eta}\|_q-\!\langle{\bm{e}_{n+k}^{n+K+1}},\bm{\eta}\rangle\!\leq 0,\;\text{for all}\; k=0,\cdots,K, \label{SOCPlinearconstraintNorm}\\
& \;\qquad\qquad\qquad -\tilde{\bm{N}}\bm{\eta} \leq \bm{0},\label{SOCPlinearconstraint2Norm}\\
& \;\qquad\qquad\qquad \|\bm{N}\bm{\eta}\|_{2} \leq 1, \label{SOCPnormconstraintNorm}
\end{align}
\label{SOCPNorm} 
\end{subequations}
where the vector inequality in 
\eqref{SOCPlinearconstraint2Norm} is elementwise.

Problem \eqref{SOCPNorm} is an approximation of \eqref{ConvexifiedSubproblemNorm}-\eqref{DefcANDgammaNorm}, which in turn is a convex relaxation of \eqref{SptFnOptimizationCollisionDetection}. Suppose the integral approximation in \eqref{SOCPNorm} incurs a local truncation error $\varepsilon>0$ w.r.t. \eqref{ConvexifiedSubproblemNorm}-\eqref{DefcANDgammaNorm}. The optimal value in \eqref{SOCPobjNorm} is therefore parameterized by $\varepsilon$. We have the following result. 

\begin{theorem}\label{thm:losslessconexificationNorm}
Let $\tilde{p}_{\varepsilon}^{*}$ be the optimal value of \eqref{SOCPNorm} where $\varepsilon>0$ is the local truncation error due to numerically approximating the integral in \eqref{ConvexifiedSubproblemNorm}. Let $\tilde{p}^{*}$ be the optimal value of \eqref{ConvexifiedSubproblemNorm}-\eqref{DefcANDgammaNorm}. Let $p^{*}$ be the optimal value obtained by solving the associated problem \eqref{SptFnOptimizationCollisionDetection}. Then,\\
(i) $\tilde{p}^{*}=\lim_{\varepsilon\downarrow 0}\tilde{p}_{\varepsilon}^{*}$.\\
(ii) $\tilde{p}^* \leq 0$.\\
(iii) $\tilde{p}^* = 0 \Rightarrow 0 \leq p^{*} \Leftrightarrow$ $\mathcal{X}_{t}^{\texttt{A}}$ and $\mathcal{X}_{t}^{\texttt{B}}$ intersect.\\
(iv) $\tilde{p}^* < 0 \Rightarrow \tilde{p}^* = p^{*} < 0 \Leftrightarrow$ $\mathcal{X}_{t}^{\texttt{A}}$ and $\mathcal{X}_{t}^{\texttt{B}}$ are disjoint.
\end{theorem}
The proof for Theorem \ref{thm:losslessconexificationNorm} is provided in Appendix \ref{AppendixProofThm:losslessconexificationNorm}. The convex problem \eqref{SOCPNorm} can be solved numerically using standard interior point algorithms\footnote{For details on computational complexity, see Supplementary Material.}.
\section{The Case when $\mathcal{U}^{\texttt{A}},\mathcal{U}^{\texttt{B}}$  are Hyperrectangles}\label{sec:HyperrectangleInputSets}
We next consider a generalized version of the $\infty$-norm bounded input uncertainties in the sense we allow hyperrectangular or box-valued input uncertainty sets of the form 
\begin{align}
\label{BoxInputSet} 
&\mathcal{U}(s):=\left[\alpha_1(s),\beta_1(s)\right]\times \left[\alpha_2(s),\beta_2(s)\right]\times \nonumber\\
&\hspace*{0.5in} \hdots \left[\alpha_m(s),\beta_m(s)\right]\subset\mathbb{R}^{m}\: \text{for all}~0\leq s \leq t,
\end{align}
where $\times$ denotes the Cartesian product. Notice that when $|\alpha(s)|=|\beta(s)|$ for all $s \in [0,t]$, the input set \eqref{BoxInputSet} represents an $\infty$-norm ball as in \eqref{inputNorm}. We will show that the intersection certification or falsification, in this case, reduces to solving $m$ decoupled second order cone programs (SOCPs) where $m$ is the number of inputs.

In this setting, 
the integrator reach set needs to account for all possible combinations of worst-cases of all the input components. Consequently, the  block diagonal system matrices as in (\ref{blkdiagAB}), makes each of the $m$ single input integrator dynamics with $r_{j}$ dimensional state subvectors for $j \in [m]$, decoupled from each other. Hence, $\mathcal{X}_t \subset \mathbb{R}^n$ is the Cartesian product of these single input integrator reach sets $\mathcal{X}_{jt}\subset \mathbb{R}^{r_{j}}$ for $j\in[m]$, i.e.,
{\small{\begin{align}
\label{CartesianProduct}
\mathcal{X}_t\!&=\mathcal{X}_{1t}\times\mathcal{X}_{2t}\times\hdots\times\mathcal{X}_{mt}\nonumber\\ 
&=\!\prod_{j=1}^m \Big\{\exp(t\bm{A}_j)\bm{x}_{0} \dotplus\!\!\int_{0}^{t}\!\!\!\bm{\xi}_j(s)\left[\alpha_j(s),\beta_j(s)\right]\differential s\Big\}, 	
\end{align}}}
where $\dotplus$ denotes the Minkowski sum, and the vectors $\bm{\xi}_{j}(s)$ are as in \eqref{xiVector}. Notice that (\ref{CartesianProduct}) may also be written as\footnote{{{In general, the Minkowski sum of a given collection of compact convex sets is not equal to their Cartesian product. However, the ``factor sets" in \eqref{CartesianProduct} belong to disjoint mutually orthogonal $r_{j}$ dimensional subspaces, $j\in[m]$, which allows writing this specific Cartesian product as a Minkowski sum.}}} a Minkowski sum $\mathcal{X}_{1t} \dotplus \hdots \dotplus \mathcal{X}_{mt}$.

Since the support function of the Minkowski sum is equal to the sum of the support functions, we have
\begin{align}
\vspace*{-0.1in}
h_{\mathcal{X}_t}(\bm{y}) = \sum_{j=1}^{m} h_{\mathcal{X}_{jt}}(\bm{y}_{j}).
\label{SptFn}    
\end{align}
where the summand support functions in the RHS of \eqref{SptFn} are given in the following theorem (proof in Appendix \ref{AppendixProofThm:SptFn}).

\begin{theorem}\label{Thm:SptFn}
The support functions of the reach sets $\mathcal{X}_{jt}\subset \mathbb{R}^{r_{j}}$ for $j\in[m]$ at time $t$ with input set $[\alpha_{j}(s),\beta_{j}(s)]$, and initial condition $\bm{x}_{j0} \in \mathbb{R}^{r_{j}}$, is\vspace{-0.1 in}
\begin{align}
 &h_{\mathcal{X}_{jt}}(\bm{y}_{j}) =\langle\bm{y}_{j}, \exp\!\left(t\bm{A}_j\right)\bm{x}_{j0}\rangle\!\nonumber\\ 
 &+\!\int_{0}^{t}\!\!\!\big[\nu_j(s) \langle \bm{y}_j,\bm{\xi}_{j}(s)\rangle+\mu_j(s) \mid\!\langle \bm{y}_j,\bm{\xi}_{j}(s) \rangle\! \mid \!\big] \differential s,
\label{SptFnIntegratorBox} 
\end{align}
where $\bm{\xi}$ is defined in \eqref{xiVector}, and
\begin{align}
\mu_{j}(s):=\!\left(\beta_{j}(s)-\alpha_{j}(s)\right)\!/2, \;\nu_{j}(s):=\!\left(\beta_{j}(s)+\alpha_{j}(s)\right)\!/2.
\label{DefMuNu}    
\end{align}
\end{theorem}

\begin{remark}\label{ReachSetOverapprox}
The formula \eqref{SptFn}-\eqref{SptFnIntegratorBox} provide upper bound for the support function of integrator reach set resulting from the same initial condition and arbitrary compact   $\mathcal{U}(s)\subset\mathbb{R}^{m}$ with $\alpha_j(s) := \underset{\bm{u}(s)\in\mathcal{U}(s)}{\min} \; u_{j}(s),~ \beta_j(s) := \underset{\bm{u}(s)\in\mathcal{U}(s)}{{\max} \; u_{j}(s)}$ for $0\leq s\leq t$. Therefore, \eqref{CartesianProduct} over-approximates the reach set with arbitrary compact $\mathcal{U}(s)$.  
\end{remark}

Define $\bm{\alpha}^{i}(s),\bm{\beta}^{i}(s),\bm{\mu}^{i}(s), \bm{\nu}^{i}(s)$ as in \eqref{BoxInputSet}, associated with respective input sets $\mathcal{U}^{i}(s)$ for the agents $i\in\{\texttt{A},\texttt{B}\}$, for all $0\leq s \leq t$.
Following Theorem \ref{Thm:SptFn}, we then obtain  $h_{\mathcal{X}_{t}^{\texttt{A}}}(\bm{y})$ and $h_{\mathcal{X}_{t}^{\texttt{B}}}(\bm{y})$. It then remains to solve \eqref{SptFnOptimizationCollisionDetection}.
\subsection{Distributed Computation}\label{subsec:DisComp}
Instead of directly substituting $h_{\mathcal{X}_{t}^{\texttt{A}}}(\bm{y})$ and $h_{\mathcal{X}_{t}^{\texttt{B}}}(\bm{y})$ in \eqref{SptFnOptimizationCollisionDetection}, we make the observation that 
{\small{\begin{align}
\label{ProductIntersectionInterchange}
\!\!\mathcal{X}_{t}^{\texttt{A}}\cap\mathcal{X}_{t}^{\texttt{B}} =\! \left(\prod_{j=1}^{m}\mathcal{X}_{jt}^{\texttt{A}}\!\right) \!\!\cap\!\! \left(\prod_{j=1}^{m}\mathcal{X}_{jt}^{\texttt{B}}\!\right) \!\!=\! \prod_{j=1}^{m} \left(\mathcal{X}_{jt}^{\texttt{A}} \cap \mathcal{X}_{jt}^{\texttt{B}}\right),
\end{align}}}
where $\prod$ denotes the Cartesian product, and $\mathcal{X}_{jt}^{\texttt{A}},\mathcal{X}_{jt}^{\texttt{B}}\subset\mathbb{R}^{r_{j}}$ are the respective $j$th single input integrator reach sets resulting from their input sets $[\alpha_{j}^{\texttt{A}}(s),\beta_{j}^{\texttt{A}}(s)]$ and $[\alpha_{j}^{\texttt{B}}(s),\beta_{j}^{\texttt{B}}(s)]$, $0\leq s \leq t$. From \eqref{ProductIntersectionInterchange}, it follows that $\mathcal{X}_{t}^{\texttt{A}}\cap\mathcal{X}_{t}^{\texttt{B}}=\varnothing$ iff there exists $j\in[m]$ such that $\mathcal{X}_{jt}^{\texttt{A}} \cap \mathcal{X}_{jt}^{\texttt{B}} = \varnothing$.

Therefore, it suffices to check whether these \emph{single input} integrator reach sets $\mathcal{X}_{jt}^{\texttt{A}},\mathcal{X}_{jt}^{\texttt{B}}$ intersect or not. Consequently, problem \eqref{SptFnOptimizationCollisionDetection} can be solved in a \emph{distributed manner}, i.e., by separately solving  
\begin{align}
\label{SptFnOptimizationCollisionDetectionDistributed} 
\underset{\bm{y}_{j}\in\mathbb{R}^{r_j},\|\bm{y}_{j}\|_{2}=1}{\min}\: h_{\mathcal{X}_{jt}^{\texttt{A}}}(\bm{y}_{j}) + h_{\mathcal{X}_{jt}^{\texttt{B}}}(-\bm{y}_{j}),
\end{align}
for all $j\in[m]$, and then checking the signs of these $m$ minimum values. In summary, $\mathcal{X}_{t}^{\texttt{A}}$ and $\mathcal{X}_{t}^{\texttt{B}}$ intersect iff \eqref{SptFnOptimizationCollisionDetectionDistributed} yields $\geq 0$ for all $j\in[m]$.

We relax the unit sphere constraint $\|\bm{y}_{j}\|_{2}=1$ to $\|\bm{y}_{j}\|_{2}\leq1$ in subproblems \eqref{SptFnOptimizationCollisionDetectionDistributed} with support functions $h_{\mathcal{X}_{jt}}$ given by \eqref{SptFnIntegratorBox}. Following the same steps as Sec. \ref{subsec:convexification}, we can rewrite subproblems \eqref{SptFnOptimizationCollisionDetectionDistributed} for each $j \in [m]$, as SOCP:
\begin{subequations}
\begin{align}
&\underset{\bm{\eta}^{{\tiny{\textup{box}}}}_{j}\in\mathbb{R}^{r_{j}+K+1}}{\min}\:\langle \bm{\kappa}_{j}^{{\tiny{\textup{box}}}}(t),\bm{\eta}_{j}^{{\tiny{\textup{box}}}}\rangle   \label{SOCPobj}\\ 
&\text{subject to} \quad \bm{M}_{j}^{{\tiny{\textup{box}}}}\bm{\eta}_{j}^{{\tiny{\textup{box}}}} \leq \bm{0}, \label{SOCPlinearconstraint}\\
& \qquad\qquad \|\bm{N}_{j}^{{\tiny{\textup{box}}}}\bm{\eta}_{j}^{{\tiny{\textup{box}}}}\|_{2} \leq 1, \label{SOCPnormconstraint}
\end{align}
\label{SOCP} 
\end{subequations}
where 
{\small{\begin{subequations}
\begin{align}
&\bm{c}_{j}^{{\tiny{\textup{box}}}}(t) :=\exp(t\bm{A}_{j})\left(\bm{x}_{j0}^{\texttt{A}} - \bm{x}_{j0}^{\texttt{B}}\right) \nonumber\\
&\qquad\qquad\qquad\qquad\qquad + \int_{0}^{t}\left(\nu_{j}^{\texttt{A}}(s)-\nu_{j}^{\texttt{B}}(s)\right)\!\bm{\xi}_{j}(s)\differential s,\label{Defc}\\
&\bm{\gamma}_{j}^{{\tiny{\textup{box}}}}(s) := \left(\mu_{j}^{\texttt{A}}(s)+\mu_{j}^{\texttt{B}}(s)\right)\bm{\xi}_{j}(s).\label{Defgamma}\\ 
&\bm{\eta}^{{\tiny{\textup{box}}}}_{j} := \begin{pmatrix}
\bm{y}_{j}\\
\bm{\theta}_{j}
\end{pmatrix}\in\mathbb{R}^{r_{j}+K+1},\label{defzj}\\
&\bm{\omega}_{j}^{{\tiny{\textup{box}}}} := \Delta s\begin{pmatrix}
1/2\\
\bm{1}_{K-1}\\
1/2
\end{pmatrix}\in\mathbb{R}^{K+1}_{>0}
,\label{defomegaj}\\ 
&\bm{\kappa}_{j}^{{\tiny{\textup{box}}}}(t) := \begin{pmatrix}
\bm{c}^{{\tiny{\textup{box}}}}_{j}(t)\\
\bm{\omega}^{{\tiny{\textup{box}}}}_{j}
\end{pmatrix}\in\mathbb{R}^{r_{j}+K+1}, \label{defellj}\\
&\bm{\Gamma}^{{\tiny{\textup{box}}}}_{j}\!:= 
\begin{pmatrix}
\bm{\gamma}_{j}^{{\tiny{\textup{box}}}\top}(s_{0})\\
-\bm{\gamma}_{j}^{{\tiny{\textup{box}}}\top}(s_{0})\\
\vdots\\
\bm{\gamma}_{j}^{{\tiny{\textup{box}}}\top}(s_{K})\\
-\bm{\gamma}_{j}^{{\tiny{\textup{box}}}\top}(s_{K})
\end{pmatrix}\in\mathbb{R}^{2(K+1)\times r_{j}},\\
&\vspace*{-0.8 in}\bm{M}_{j}^{{\tiny{\textup{box}}}} := \begin{pmatrix}
\begin{array}{c|c}
\bm{\Gamma}^{{\tiny{\textup{box}}}}_{j} & -\bm{I}_{K+1}\otimes\bm{1}_{2} \\
\hline
\bm{0}_{(K+1)\times r_{j}} & -\bm{I}_{K+1}
\end{array} 
\end{pmatrix}\in\mathbb{R}^{3(K+1)\times(r_{j}+K+1)}, \label{defMj}\\
&\bm{N}^{{\tiny{\textup{box}}}}_{j} := \!\begin{pmatrix}
\begin{array}{c|c}
\bm{I}_{r_{j}} & \bm{0}_{r_{j}\times(K+1)}
\end{array} 
\end{pmatrix}\in\mathbb{R}^{r_{j}\times(r_{j}+K+1)},\label{defNj} 
\end{align}
\label{defSOCPparam}
\end{subequations}}}
and the symbol $\otimes$ denotes the Kronecker product. 
 
We observe that \eqref{SOCPnormconstraint} results from the convexification of the nonconvex constraint $\|\bm{y}_{j}\|_{2}=1$ in \eqref{SptFnOptimizationCollisionDetectionDistributed} for each $j\in[m]$. As in Sec. \ref{subsec:convexification},  convexification of \eqref{SptFnOptimizationCollisionDetectionDistributed} turns out to be lossless. In the small $\Delta s$ limit, the local truncation error $\varepsilon$ in approximating the integral w.r.t. $t$, goes to zero and numerically solving \eqref{SOCP} allows us to certify $\mathcal{X}_{jt}^{\texttt{A}}\cap\mathcal{X}_{jt}^{\texttt{B}} \neq (=) \varnothing$. We summarize this in the following statement whose proof follows the steps as in Appendix \ref{AppendixProofThm:losslessconexificationNorm}, and is omitted.
\begin{theorem}\label{thm:losslessconexificationBox}
For $j\in[m]$, let $p_{j}^*$ be the optimal value of problem \eqref{SptFnOptimizationCollisionDetectionDistributed}, and let $\tilde{p}_{j}^*$ be the optimal value for its convexification by replacing $\|\bm{y}_{j}\|_{2}=1$ with $\|\bm{y}_{j}\|_{2}\leq 1$. Let $\tilde{p}_{j,\varepsilon}^{*}$ be the optimal value of \eqref{SOCP}-\eqref{defSOCPparam} where $\varepsilon>0$ is the local truncation error due to numerically approximating the integral. Then, the following holds:\\
(i) $\tilde{p}_{j}^* = \lim_{\varepsilon\downarrow 0}\tilde{p}_{j,\varepsilon}^*$.\\
(ii) $\tilde{p}_{j}^* \leq 0$.\\
(iii) $\tilde{p}_{j}^* = 0 \Rightarrow 0 \leq p_{j}^{*} \Leftrightarrow$ $\mathcal{X}_{jt}^{\texttt{A}}$ and $\mathcal{X}_{jt}^{\texttt{B}}$ intersect.\\
(iv) $\tilde{p}_{j}^* < 0 \Rightarrow \tilde{p}_{j}^* = p_{j}^{*} < 0 \Leftrightarrow$ $\mathcal{X}_{jt}^{\texttt{A}}$ and $\mathcal{X}_{jt}^{\texttt{B}}$ are disjoint.
\end{theorem}

Since \eqref{ProductIntersectionInterchange} tells us $\mathcal{X}_{t}^{\texttt{A}}\cap\mathcal{X}_{t}^{\texttt{B}}=\varnothing$ iff there exists $j\in[m]$ such that $\mathcal{X}_{jt}^{\texttt{A}} \cap \mathcal{X}_{jt}^{\texttt{B}} = \varnothing$, therefore the distributed computation of \eqref{SOCP} allows certification or falsification of integrator reach sets subject to box-valued input uncertainties. 

In the following, we provide a numerical example to illustrate our results.

\begin{figure}[t]
    \centering
    \includegraphics[width=0.9\linewidth]{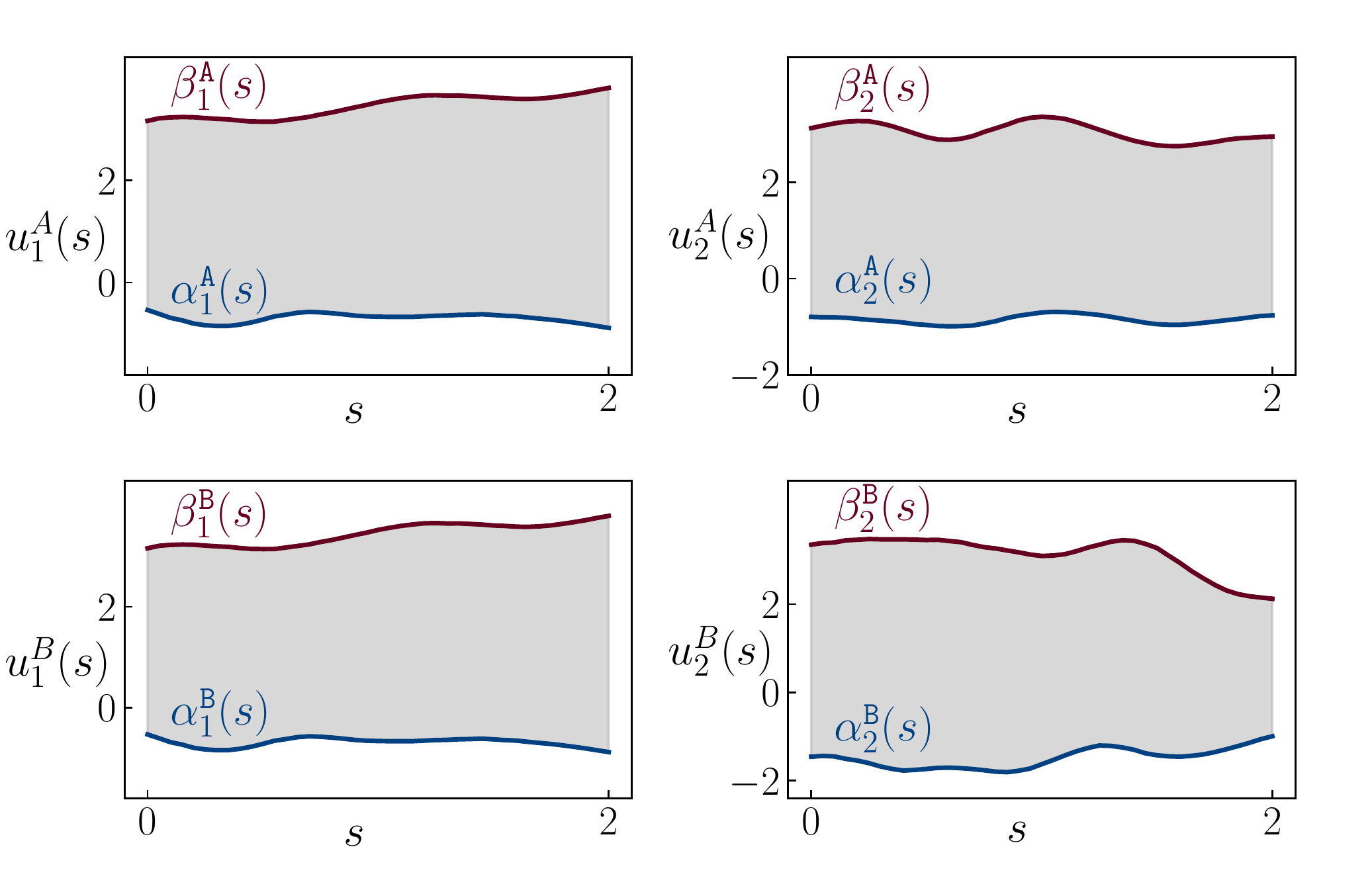}
    \caption{{\small{The input uncertainties for agent $\texttt{A}$ with $u_j^\texttt{A}(s)\in[\alpha_{j}^{\texttt{A}}(s),\beta_{j}^{\texttt{A}}(s)]$, and agent $\texttt{B}$ with $u_j^\texttt{B}(s) \in [\alpha_{j}^{\texttt{B}}(s),\beta_{j}^{\texttt{B}}(s)]$, where $j\in\{1,2\}$ for the example in Sec. \ref{example:CollisionDetection}. The controlled dynamics of the agents are given by \eqref{Example3block2block}. At each $s \in [0,2]$, ${\alpha}^i_j(s)$ and ${\beta}^i_j(s)$ respectively represent the  minimum and maximum of the $j$th coordinate of the input set $\mathcal{U}^{i}(t)$ for $i\in\{\texttt{A}, \texttt{B}\}$.}}}
\vspace*{-0.23in}
\label{Example1InputSets}
\end{figure}
\begin{figure}[th]
    \centering
    \includegraphics[width=0.89\linewidth]{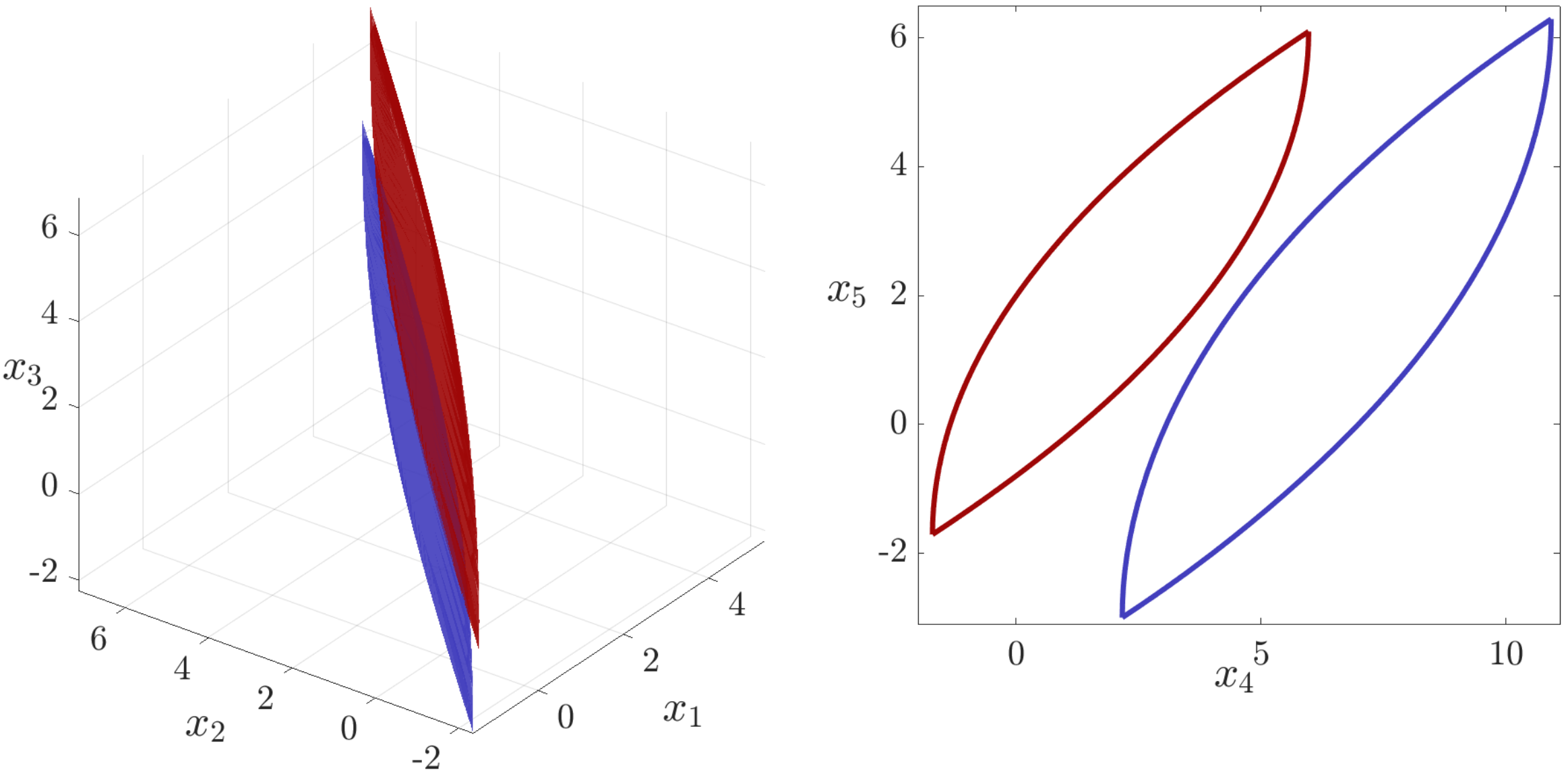}
    \caption[Caption for LOF]{{\small{Intersection detection between two integrator reach sets $\mathcal{X}_{t}^{\texttt{A}}$ (\emph{red}) and $\mathcal{X}_{t}^{\texttt{B}}$ (\emph{blue}) described in Sec. \ref{example:CollisionDetection} with relative degree vectors $\bm{r}^{\texttt{A}}=\bm{r}^{\texttt{B}}=(3,2)^{\top}$. The agents start from initial conditions $\bm{x}_0^{\texttt{A}}=(0.5,\bm{0}_{{\small{1\times4}}})^{\top}$ and $\bm{x}_0^{\texttt{B}}=(\bm{0}_{{\small{1\times3}
    }}, 5,0)^{\top}$. In this case, $\mathcal{X}_{t}^{\texttt{A}}=\mathcal{X}_{1t}^{\texttt{A}}\times\mathcal{X}_{2t}^{\texttt{A}}, \mathcal{X}_{t}^{\texttt{B}}=\mathcal{X}_{1t}^{\texttt{B}}\times\mathcal{X}_{2t}^{\texttt{B}}\subset\mathbb{R}^{5}$. \emph{Left}: Intersection between the sets $\mathcal{X}_{1t}^{\texttt{A}}$ (\emph{red}) and $\mathcal{X}_{1t}^{\texttt{B}}$ (\emph{blue}). \emph{Right}: no intersection between the sets $\mathcal{X}_{2t}^{\texttt{A}}$ (\emph{red}) and $\mathcal{X}_{2t}^{\texttt{B}}$ (\emph{blue}). These plots are made by generalizing the parametric boundary formula of the integrator reach sets in \cite[proposition 1]{haddad2021curious} to account for time-varying set-valued uncertainties.}}}
\vspace*{-0.23in}
\label{Example1Intersection}
\end{figure}

\subsection{Numerical Example}\label{example:CollisionDetection}
Let us consider two integrator agents \texttt{A} and \texttt{B} with relative degrees $\bm{r}^{\texttt{A}}=\bm{r}^{\texttt{B}}=(3,2)^{\top}$ as in \eqref{Example3block2block}, starting from respective initial conditions $$\bm{x}_0^{\texttt{A}}=(0.5,\bm{0}_{{\small{1\times4}}})^{\top},\quad \bm{x}_0^{\texttt{B}}=(\bm{0}_{{\small{1\times3}}}, 5,0)^{\top}.$$ 
Given two box-valued input sets $\mathcal{U}^{\texttt{A}}(s)=[\alpha_{1}^{\texttt{A}}(s),\beta_{1}^{\texttt{A}}(s)]\times [\alpha_{2}^{\texttt{A}}(s),\beta_{2}^{\texttt{A}}(s)]$ and $\mathcal{U}^{\texttt{B}}(s)=[\alpha_{1}^{\texttt{B}}(s),\beta_{1}^{\texttt{B}}(s)]\times [\alpha_{2}^{\texttt{B}}(s),\beta_{2}^{\texttt{B}}(s)]$ shown in Fig. \ref{Example1InputSets}, we want to check if there is an intersection at time $t=2$ between the reach sets of agent $\texttt{A}$, denoted as $\mathcal{X}_{t}^{\texttt{A}}$, and that of agent $\texttt{B}$, denoted as $\mathcal{X}_{t}^{\texttt{B}}$. 

As explained in Sec. \ref{subsec:DisComp}, it suffices to check if there will be an intersection between each corresponding $j$th single input reach sets $\mathcal{X}_{jt}^{\texttt{A}}$ and $\mathcal{X}_{jt}^{\texttt{B}}$, for $j\in\{1,2\}$. 

Using (\ref{defSOCPparam}), we construct the matrices $\bm{M}_j^{i,{\tiny{\textup{box}}}}, \bm{N}_j^{i,{\tiny{\textup{box}}}}$ and the vector $\bm{\kappa}_j^{i{\tiny{\textup{box}}}}$ for $i\in\{\texttt{A},\texttt{B}\}$. Then, we solve the optimization problems (\ref{SOCP}) for each $j\in\{1,2\}$ via MATLAB CVX toolbox \cite{cvx,gb08} with $\Delta s=0.05$. The runtimes are $0.38$ s and $0.37$ s for $j=1$ and $j=2$, respectively. The corresponding optimal values are 
\begin{align}
(\tilde{p}_{1}^*,\tilde{p}_{2}^*)=(0,-0.54).
\label{tildestarvalues}	
\end{align}
These optimal values imply $\mathcal{X}_{1t}^{\texttt{A}} \cap \mathcal{X}_{1t}^{\texttt{B}} \ne \varnothing$, and $\mathcal{X}_{2t}^{\texttt{A}} \cap \mathcal{X}_{2t}^{\texttt{B}} = \varnothing$. Therefore, we conclude: $\mathcal{X}_{t}^{\texttt{A}} \cap \mathcal{X}_{t}^{\texttt{B}} = \varnothing$. The same pair \eqref{tildestarvalues} was obtained for smaller step-sizes:
\begin{itemize}
\item $\Delta s=0.01$, runtimes are 0.80 s for both $j = 1$ and $2$.
\item $\Delta s=0.005 $, runtimes are 1.24 s and 1.25 s for $j = 1$ and $j = 2$, respectively.
\end{itemize}

In agreement with Theorem \ref{thm:losslessconexificationBox}, for $\tilde{p}_{1}^*=0$, CVX returns $\|\bm{N}_{1}\bm{\eta}_{1}\|_{2}=0$, while for $\tilde{p}_{2}^*<0$, it returns $\|\bm{N}_{2}\bm{\eta}_{2}\|_{2}=1$ demonstrating that 
the convexification (\ref{SOCPnormconstraint}) is lossless. This scenario is illustrated in Fig. \ref{Example1Intersection}. 





\section{Conclusions}\label{sec:conclusions}
This work presents a variational formulation for certifying or falsifying intersection of the reach sets of integrator agents subject to set-valued input uncertainties. The proposed nonconvex formulation is shown to enjoy lossless convexification for time-varying norm bounded as well as generic hyperrectangular input uncertainties, thus being amenable to convex programming for tractable computation. A numerical example is provided to illustrate the results.

\appendix
\subsection{{Proof of Proposition \ref{proSptFn}}}\label{AppendixProofPro:proSptFn}
Support function is distributive over sum, so from \eqref{IntegratorODE} and \eqref{DefSptFn}, we obtain 
\begin{align}
\label{hRjgeneral}
h_{\mathcal{X}_{t}}\left(\bm{y}\right) = &\langle\bm{y},\exp\left(t\bm{A}\right)\bm{x}_{0}\rangle+ h_{\int_{0}^{t}\exp(s\bm{A})\bm{B}\:\mathcal{U}(s)\:\differential s}(\bm{y}).
\end{align}
Using \cite[Proposition 1]{haddad2020convex}, we then have
\begin{align}
\label{hrjprop}
&h_{\int_{0}^{t}\exp(s\bm{A})\bm{B}\:\mathcal{U}(s)\:\differential s}(\bm{y}) = \int_{0}^{t} \!\!h_{\exp(s\bm{A})\bm{B}\: \mathcal{U}(s)}\:(\bm{y})\;\differential{s}\nonumber\\ 
&\qquad\qquad =\int_{0}^{t}\!\!\underset{\bm{u}(s) \in \mathcal{U}(s) }{\sup}\!\langle \bm{y}_{},\exp(s\bm{A})\bm{B}\bm{u}(s)\rangle~ \differential s. 
\end{align}
Combining \eqref{hrjprop} with the structures of the state and input matrices in \eqref{IntegratorODE}-\eqref{blkdiagAB}, allows us to rewrite \eqref{hRjgeneral} as \eqref{SptFnIntegratorpGen}. Specifically, the result follows from the fact that the state transition matrix $\exp(s\bm{A})=\blkdiag\left(\exp(s\bm{A}_{1}),\hdots,\exp(s\bm{A}_{m})\right)$ wherein each diagonal block is upper triangular with entries
$$\left(\exp(s\bm{A}_{j})\right)_{ab} = \begin{cases}
\dfrac{s^{a-b}}{(a-b)!}&\text{for}\;a\leq b,\\
0 &\text{otherwise},
\end{cases}$$
for all $a,b=1,\hdots,r_{j}$, for each $j\in[m]$.
\hfill\qed
\subsection{{{Proof of Theorem \ref{thm:Normspf}}}}\label{AppendixProofThm:Normspf}
Let $\bm{z}(s):=(\exp(s\bm{A})\bm{B})^{\top}\bm{y}$ for all $0\leq s \leq t$. Then, the integrand of the RHS of \eqref{hrjprop} for time-varying $p$-norm ball input set $\mathcal{U}(s)$ as in \eqref{inputNorm}, becomes
\begin{align}
&\underset{\bm{u}(s)\in\mathbb{R}^{m}}{\sup} \quad \:\langle \bm{z}(s),\bm{u}(s)\rangle \nonumber\\ 
&\text{subject to} \quad \|\bm{u}(s)\|_p \leq \ell(s).
\label{LP} 
\end{align}
Recall that $\|\cdot\|_{q}$ is the dual norm of $\|\cdot\|_{p}$ for $\frac{1}{p}+\frac{1}{q}=1$. From the definition of dual norm, we have $$\langle \bm{z}(s),\bm{u}(s)\rangle \leq \|\bm{z}(s)\|_{q} \|\bm{u}(s)\|_{p} \overset{\eqref{LP}}{\leq} \ell(s)\|\bm{z}(s)\|_{q}$$
with equality resulting in the supremum in \eqref{LP}. Substituting this supremum in \eqref{hrjprop} and then using \eqref{hRjgeneral}, we get \eqref{SptFnIntegratorpNorm}. \hfill\qed

\subsection{{{Proof of Theorem \ref{thm:losslessconexificationNorm}}}}\label{AppendixProofThm:losslessconexificationNorm}
(i) Follows immediately from \eqref{TrapzApprox}, \eqref{defSOCPparamNorm} and \eqref{ConvexifiedSubproblemNorm}.

(ii) Since $\bm{0}_{n\times 1}$ is in the feasible set of \eqref{ConvexifiedSubproblemNorm} and makes the objective equal to zero, we have $\tilde{p}^*\leq 0$.

(iii)  Since $\tilde{p}^*$ is the optimal value of the convex relaxation of a nonconvex problem with optimal value $p^{*}$, we must have $\tilde{p}^* \leq p^{*}$. Thus, $\tilde{p}^* = 0 \Rightarrow 0 \leq p^{*}$. That $0\leq p^{*}$ is equivalent to $\mathcal{X}_{t}^{\texttt{A}}\cap\mathcal{X}_{t}^{\texttt{B}} \neq\varnothing$, was explained before in Sec. \ref{sec:formulation}.

(iv) We now show that when $\tilde{p}^* < 0$, the convexification is in fact lossless, i.e., $\tilde{p}^* = p^{*}$. Denote the $\arg\min$ for the convex problem \eqref{ConvexifiedSubproblemNorm} as $\bm{y}^{\min}$. It suffices to prove that $\|\bm{y}^{\min}\|_2 = 1$. To this end, suppose if possible, that $\|\bm{y}^{\min}\|_{2} =: \delta < 1$, i.e., $0<\delta < 1$. Now let $\widetilde{\bm{y}} := \bm{y}^{\min}/\delta$, which is clearly feasible w.r.t. \eqref{ConvexifiedSubproblemNorm}. However,
\[f_{0}\left(\widetilde{\bm{y}}\right) = \underbrace{\frac{1}{\delta}}_{>1} \underbrace{\vphantom{\frac{1}{\varepsilon}} f_0(\bm{y}^{\min})}_{<0} < f_0(\bm{y}^{\min}),\]
contradicting the supposition that the $\arg\min$ is strictly within the unit ball in $\mathbb{R}^{n}$. Therefore, if $\tilde{p}^{*}<0$ then $\tilde{p}^{*} = {p}^{*}$. That $p^{*}<0\Leftrightarrow \mathcal{X}_{t}^{\texttt{A}}\cap\mathcal{X}_{t}^{\texttt{B}} =\varnothing$, was explained in Sec. \ref{sec:formulation}.\hfill\qed

\subsection{{Proof of Theorem \ref{Thm:SptFn}}}\label{AppendixProofThm:SptFn}
Specializing \eqref{SptFnIntegratorpGen} for the box-valued input set \eqref{BoxInputSet}, for all $j\in[m]$, we get 
{\small{\begin{align}
h_{\mathcal{X}_{jt}}\left(\bm{y}_{j}\right) &= 
 \langle\bm{y}_{j},\exp\left(s\bm{A}_{j}\right)\bm{x}_{j0}\rangle \nonumber\\ 
&+ \int^t_0 \underset{u_j(s)\in\mathcal[\alpha_j(s),\beta_j(s)] }{\sup} \langle\bm{y}_{j},\bm{\xi}_{j}(s)\rangle\: u_{j}(s)\: \differential{s}.
\label{SptFnbox}
\end{align}}}	

The optimizer of the integrand in the RHS of \eqref{SptFnbox}, $u_{j}^{\text{opt}}$, can be written in terms of the Heaviside unit step function $H(\cdot)$ as 
\begin{align}
u_{j}^{\text{opt}} &= \alpha_j + (\beta_j - \alpha_j) H(\langle\bm{y}_{j},\bm{\xi}_{j}\rangle)\nonumber\\
&= \alpha_j + (\beta_j - \alpha_j) \times \frac{1}{2}\left(1 + {\rm{sgn}}\left(\langle\bm{y}_{j},\bm{\xi}_{j}\rangle\right)\right),
\label{minimizerspt}
\end{align}
where ${\rm{sgn}}(\cdot)$ denotes the sign function. Substituting the optimizer \eqref{minimizerspt} into \eqref{SptFnbox} yields \eqref{SptFnIntegratorBox}.\hfill\qed

\bibliographystyle{IEEEtran}
\bibliography{References.bib}

\end{document}


\maketitle
\pagenumbering{arabic}



The purpose of this Supplementary Material is to provide some details supporting the main text.

\section*{Related Literature}\label{sec:Literature}

Motivated by safety-critical applications, there exists a vast literature on reach set computation, see e.g., \cite{althoff2021set,chutinan1999verification,villanueva2015unified,le2010reachability,pecsvaradi1971reachable}. In particular, one may be interested to detect the intersection between two reach sets in specific contexts such as collision avoidance in automated driving \cite{gruber2018anytime,1253216}. 

Existing algorithms for detecting the intersection of sets--a problem that also arises in computer graphics \cite{466717,cpg}--use either the distance formulation or the support function formulation (see e.g., \cite{lin2002distance,iwata2015computing}) mentioned in the Sec. III of the main text. However, it is not \emph{a priori} clear which of the two formulations is computationally more benign for specific sets. One of the contributions of our work is to show that the integrator reach sets being zonoids \cite{haddad2021curious}, the support function based variational formulation is computationally tractable for certifying or falsifying their intersections.

\section*{Intersection of Reach Sets of Static State Feedback Linearizable Systems}\label{sec:StaticStateFeedback}
Full state static feedback linearizable systems of the form $\dot{\bm{z}}=\bm{f}(\bm{z},\bm{v})$ with $n\times 1$ state vector $\bm{z}$ and $m\times 1$ input vector $\bm{v}$, comprise a class of controlled dynamics that can be brought exactly in the integrator a.k.a. Brunovsky canonical form via diffeomorphic change of state co-ordinates. In other words, for this class of systems, there exist local diffeomorphism\footnote{A differentiable bijective mapping with differentiable inverse.} $$\bm{\tau}:\mathcal{Z}\subseteq \mathbb{R}^{n} \mapsto \mathcal{X}\subseteq \mathbb{R}^{n}$$ 
such that the new state $\bm{x}=\bm{\tau}(\bm{z})$ follows (see e.g., \cite[p. 243-248]{isidori1989nonlinear}) state dynamics governed by the Brunovsky normal form given by equation (2) in the main text. 

The original state-control pair $(\bm{z},\bm{v})\in\mathcal{Z}\times\mathcal{V}\subseteq\mathbb{R}^{n}\times\mathbb{R}^{m}$ gets mapped to the new state-control pair $(\bm{x},\bm{u})\in\mathcal{X}\times\mathcal{U}\subseteq\mathbb{R}^{n}\times\mathbb{R}^{m}$. The new control $\bm{u}$ is guaranteed to be a homeomorphism\footnote{A continuous bijective mapping with continuous inverse.} of the original control $\bm{v}$, and hence for all $0\leq s\leq t$, the input sets $\mathcal{U}(s)$ are compact (for any $s$) provided the input sets $\mathcal{V}(s)$ are. As we mentioned in the main text's Sec. I, fourth paragraph, since the reach set $\mathcal{X}$ generated by (1) remains invariant under the the closure of convexification of $\mathcal{U}(s)$, the results of this paper can be used to certify or falsify the intersection of the reach sets of two full state static feedback linearizable agents with identical dynamics. Let us explain this in some detail since this may help understand the scope and broader motivation for the work in the main text.

Consider two static state feedback linearizable systems labeled as \texttt{A} and \texttt{B}, each with $n$ states, $m$ inputs and relative degree $\bm{r}=(r_1,\hdots,r_m)^{\top}$, having identical dynamics 
\begin{align}
\dot{\bm{z}}^{\texttt{A}} = \bm{f}\left(\bm{z}^{\texttt{A}},\bm{v}^{\texttt{A}}\right), \quad  \dot{\bm{z}}^{\texttt{B}} = \bm{f}\left(\bm{z}^{\texttt{B}},\bm{v}^{\texttt{B}}\right).    
\label{DiffFlatDyn}    
\end{align}
Suppose we want to check whether their reach sets $\mathcal{Z}_{t}^{\texttt{A}}, \mathcal{Z}_{t}^{\texttt{B}}\subset\mathbb{R}^{n}$ at a given time $t$, resulting from respective initial conditions $\bm{z}_{0}^{\texttt{A}},\bm{z}_{0}^{\texttt{B}}\in\mathbb{R}^{n}$, and respective compact input set-valued uncertainties 
$$\mathcal{V}^{\texttt{A}}(s),\mathcal{V}^{\texttt{B}}(s)\subset\mathbb{R}^{m}, \quad 0\leq s\leq t,$$ 
intersect or not. Determining so is cumbersome because $\mathcal{Z}_{t}^{\texttt{A}}, \mathcal{Z}_{t}^{\texttt{B}}$ are compact but nonconvex in general. However, checking if $\mathcal{Z}_{t}^{\texttt{A}}\cap\mathcal{Z}_{t}^{\texttt{B}} = \varnothing$ or not, is  \emph{equivalent} to checking if the corresponding $\mathcal{X}_{t}^{\texttt{A}}, \mathcal{X}_{t}^{\texttt{B}}\subset\mathbb{R}^{n}$ in the Brunovsky normal coordinates, intersect or not. This is because 
$$\mathcal{X}_{t}^{\texttt{A}} = \bm{\tau}\left(\mathcal{Z}_{t}^{\texttt{A}}\right), \quad \mathcal{X}_{t}^{\texttt{B}} = \bm{\tau}\left(\mathcal{Z}_{t}^{\texttt{B}}\right),$$
where $\bm{\tau}$ is the known diffeomorphism associated with $\bm{f}$, and $\mathcal{X}_{t}^{\texttt{A}}, \mathcal{X}_{t}^{\texttt{B}}$ are integrator reach sets associated with initial conditions $\bm{x}_{0}^{\texttt{A}}:=\bm{\tau}\left(\bm{z}_{0}^{\texttt{A}}\right), \bm{x}_{0}^{\texttt{B}}:=\bm{\tau}\left(\bm{z}_{0}^{\texttt{B}}\right)$, and transformed compact input set-valued trajectories $\mathcal{U}^{\texttt{A}}(s), \mathcal{U}^{\texttt{B}}(s)\subset\mathbb{R}^{m}$ corresponding to the original input set-valued trajectories $\mathcal{V}^{\texttt{A}}(s), \mathcal{V}^{\texttt{B}}(s)\subset\mathbb{R}^{m}$, $0\leq s \leq t$.

Observe that 
\begin{align}
\mathcal{Z}_{t}^{\texttt{A}}\cap\mathcal{Z}_{t}^{\texttt{B}} \neq (=) \varnothing \:\Leftrightarrow\: \bm{\tau}\left(\mathcal{Z}_{t}^{\texttt{A}}\cap\mathcal{Z}_{t}^{\texttt{B}}\right) \neq (=) \varnothing
\label{tauofintersection}	
\end{align} 
since $\bm{\tau}$ is a diffeomorphism. On the other hand, being a diffeomorphism, $\bm{\tau}$ is injective, and we have (see e.g., \cite[p. 388]{lee2010introduction}): 
\begin{align}
\bm{\tau}\left(\mathcal{Z}_{t}^{\texttt{A}}\cap\mathcal{Z}_{t}^{\texttt{B}}\right) = \bm{\tau}\left(\mathcal{Z}_{t}^{\texttt{A}}\right)\cap\bm{\tau}\left(\mathcal{Z}_{t}^{\texttt{B}}\right) =  \mathcal{X}_{t}^{\texttt{A}}\cap\mathcal{X}_{t}^{\texttt{B}}.
\label{intersectionoftau}	
\end{align}
Therefore, combining \eqref{tauofintersection} and \eqref{intersectionoftau}, we obtain
\begin{align}
 \mathcal{X}_{t}^{\texttt{A}}\cap\mathcal{X}_{t}^{\texttt{B}} \neq (=) \varnothing
  \:\Leftrightarrow\:
 \mathcal{Z}_{t}^{\texttt{A}}\cap\mathcal{Z}_{t}^{\texttt{B}} \neq (=) \varnothing,
\label{sta}
\end{align}
thereby showing the equivalence for the reach set intersection detection or falsification among the integrator and original state co-ordinates.


\section*{Computational Scalability}\label{sec:CompScalability}
Here we discuss the scalability of the proposed variational computation with respect to the dimensionality of the state space.

\subsection*{When $\mathcal{U}^{\texttt{A}},\mathcal{U}^{\texttt{B}}$  are Norm Balls}
In Sec. IV of the main text, we consider the input uncertainties to be $\ell_{p}$ norm bounded for given $p\in(0,\infty]$ (see Remark 2 in the main text). We then derive a convex problem (equation (23) in main text) with linear objective in $N_d:=n+K+1$ decision variables, $K+1$ linear constraints, and $K+2$ nonlinear constraints. Out of the nonlinear constraints, one is a second-order cone constraints of dimension $K+1$ and the remaining $K+1$ are the $q$ norm constraints, where $q$ is the H\"{o}lder conjugate of $p$.

When $p$ and hence is $q$ is rational, then we can convert  \cite[Sec. 2.3(g)]{alizadeh2003second} each of the $q$ norm constraints (equations (23b) in the main text) to a $(m+1)$th dimensional second order cone constraint and a linear constraint, and recast this convex problem as a second order cone program (SOCP). Furthermore, the $K+1$ linear constraints (equations (23b) in the main text) can be viewed as second-order cone constraints of dimension 1. Thus we have total $N_c:=3(K+1)+1$ constraints:
\begin{itemize}
\item $1$ second-order cone constraint of dimension $K+1$.
\item $K+1$ second-order cone constraints of dimension $m+1$.
\item $2(K+1)$ linear/second-order cone constraints of dimension $1$.
\end{itemize}
Let $D:=\sum_{i=1}^{N_c} n_i=(K+1)(m+4)$, where $n_i$ is the dimension of the $i$th conic constraint. 

We now comment on the complexity of solving the SOCP via standard primal-dual interior point algorithms. Following  \cite{lobo1998applications}, the number of iterations to decrease the duality gap of the primal-dual interior point algorithm to a constant fraction of itself is upper bounded by $\mathcal{O}\left(\sqrt {N_c}\right)$. The computational complexity per iteration is $$\mathcal{O}\left(N_d^2D\right)=\mathcal{O}\left((n+K+1)^2(K+1)(m+4) \right).$$
Thus, for $m,K$ fixed, the complexity is quadratic in the dimension of the state space $n$.

\subsection*{When $\mathcal{U}^{\texttt{A}},\mathcal{U}^{\texttt{B}}$  are Hyperrectangles}\

Recall that the relative degree vector $(r_1,\hdots,r_m)^{\top}$ satisfies $r_1 + \hdots + r_m = n$, where $n$ is the number of states. In Sec. V of the main text, we consider the input uncertainties to be generalized hyperrectangles, and derive $m$ distributed SOCPs where $m$ is the number of inputs. 

For $j\in[m]$, the $j$th such SOCP has $N_{dj}:=r_{j}+K+1$ decision variables, and $N_{cj}:=3K+4$ constraints: 
\begin{itemize}
\item $3K+3$ linear/second-order cone constraints of dimension $1$.
\item  $1$ second order cone constraint of dimension $K+1$. 
\end{itemize}
As before, let $n_i$ denote the dimension of the $i$th conic constraint. Then, for the $j$th SOCP, $j\in[m]$, the computational complexity per iteration is $$\mathcal{O}\left(N_{dj}^2\sum_{i=1}^{N_{cj}}n_i\right)=\mathcal{O}\left((r_j+K+1)^2(4K+4) \right).$$
For $K$ fixed, the per-iteration complexity for the $j$th distributed SOCP is quadratic in the relative degree $r_{j}$, $j\in[m]$.



\bibliographystyle{IEEEtran}
\bibliography{SuppReferences.bib}